\documentclass[10pt,reprint]{revtex4-1}

\usepackage{graphicx,epsfig,amssymb,amsbsy,amsmath}
\newcommand{\bb}[1]{\left({#1}\right)}					
\newcommand{\cc}[1]{\left\{#1\right\}}					

\newcommand{\unitstep}[1]{{\rm H}(#1)}

\newcommand{\fref}[1]{Fig.~\ref{#1}}
\newcommand{\Fref}[1]{Fig.~\ref{#1}}

\newcommand{\Frefs}[1]{Figs.~\ref{#1}}
\newcommand{\eref}[1]{(\ref{#1})}

\begin{document}
\title{Non-determinism in the limit of nonsmooth dynamics}
\author{Mike R Jeffrey}
\affiliation{Department of Engineering Mathematics, University of Bristol, Queen's Building, Bristol BS8 1TR, UK}
%
\date{\today}

\begin{abstract}
Discontinuous time derivatives are used to model threshold-dependent switching in such diverse applications as dry friction, electronic control, and biological growth. In a continuous flow, a discontinuous derivative can generate multiple outcomes from a single initial state. Here we show that well defined solution sets exist for flows that become multi-valued due to grazing a discontinuity. Loss of determinism is used to quantify dynamics in the limit of infinite sensitivity to initial conditions, then applied to the dynamics of a superconducting resonator and a negatively damped oscillator. 
\end{abstract}

\maketitle


The dynamics of physical systems affected by sudden systemic changes can be modeled by  isolated discontinuities in systems of ordinary differential equations. These are increasingly applied to study the abrupt activation of processes for biological or physical control, for example relays in electronic circuits, mitosis of living cells, and opening of ion channels (see for example \cite{b99,f88,bs08,bc08}). 

Attempting to solve for dynamics at a discontinuity can lead to ambiguities in forward or backward time. The latter arise when solutions stick to the locus of discontinuity \cite{f88}, and have non-unique histories but well determined forward evolution, with physical interpretations such as frictional sticking and chattering. Ambiguities in forward time, on the other hand, imply non-deterministic evolution, and their application is less well understood. 

The purpose of this letter is to formalize a phenomenon currently emerging in dynamical systems, whereby a localized loss of determinism leads to extreme forms of bifurcations and chaos. The salient feature, presented here, is that it is generically possible for an initially well determined flow to evolve onto a point where it becomes multi-valued. The multiple outcomes form a set-valued flow that is itself well determined. We characterise these as {\it explosions} of allowed states of the system, which imply infinite sensitivity to initial conditions in the limit of infinite rate of change in the flow velocity. Such explosions may be robust, or may appear fleetingly as a parameter is varied, as determined by the causative singularity. The latter provide an explanation for bifurcation-like transitions classified in \cite{jh09}, and the former generalize a non-deterministic form of chaos discovered in \cite{cj10}. 

Two applications are presented for illustration: an experimentally derived model of a sensor whose dynamics changes abruptly between normal and super conducting temperature ranges, and an abstract model of a particle subject to direction dependent forcing and damping. 

A prototype for such systems is given by considering an $n$-dimensional state  $\bf x$, evolving according to a set of differential equations with a discontinuity, given by
\begin{equation}\label{pws}
\frac{d\;}{dt}{\bf x}={\bf f}({\bf x})=\left\{\begin{array}{lll}{\bf f}^+({\bf x})&\mbox{if}&\sigma({\bf x})>0,\\{\bf f}^-({\bf x})&\mbox{if}&\sigma({\bf x})<0,\end{array}\right.
\end{equation}
where ${\bf f}^\pm$ and $\sigma$ are smooth vector and scalar functions respectively, and the discontinuity takes place when $\sigma=0$.  
It is important to note, throughout this article, that a {\it tangency} will refer strictly to a quadratic tangency, so that if $\bf f^\pm$ lies tangent to the discontinuity surface then ${\bf f}^\pm\cdot\partial_{\bf x}\sigma=\sigma=0$, but $({\bf f}^\pm\cdot\partial_{\bf x})^2\sigma\neq0$. 

\begin{figure}[h!]\begin{center}\includegraphics[width=0.34\textwidth]{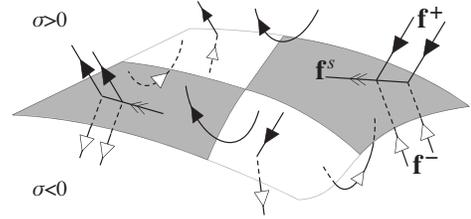}
\end{center}\vspace{-0.35cm}\caption{\sf A vector field jumps between ${\bf f}^+$ above (filled arrows), ${\bf f}^-$ below (unfilled arrows), and ${\bf f}^{s}$ inside (double arrows), a discontinuity surface. Sticking occurs in regions where ${\bf f}^\pm$ are in opposition (shaded), and crossing occurs otherwise (unshaded), with ${\bf f}^\pm$ tangent to the surface at the boundaries.}\label{fig:big2fold}\end{figure}

In general, the flow in the regions $\sigma>0$ or $\sigma<0$ is given by the first integral of \eref{pws}. Although \eref{pws} is undefined at $\sigma=0$, continuous flow solutions can be found even at the discontinuity, and the common interpretations all arrive at the same result. Near the discontinuity surface we assume that the state $\bf x$ lies in $\sigma>0$ for a fraction $\lambda$ of a small time interval $\delta t$, and lies in $\sigma<0$ for the remaining time interval $(1-\lambda)\delta t$. Then the change in  $\bf x$ over $\delta t$ is $\delta {\bf x}=\lambda\delta t {\bf f}^++(1-\lambda)\delta t {\bf f}^-$. Letting $\delta t\rightarrow0$ and allowing $\lambda$ to take any real value, we therefore augment \eref{pws} with the differential inclusion
\vspace{-0.2cm}\begin{equation}\label{inclusion}
\frac{d\;}{dt}{\bf x}\in\lambda {\bf f}^++(1-\lambda){\bf f}^-,\quad\lambda\in\mathbb R\;\;\mbox{if}\;\;\sigma=0.\vspace{-0.1cm}
\end{equation}
If the flow passes through the discontinuity, this multi-valued equation applies only at the instant when $\sigma=0$, and allows orbits in $\sigma>0$ and $\sigma<0$ to be concatenated. Thus a unique flow {\it crossing} the discontinuity is formed, illustrated at the unshaded regions in \fref{fig:big2fold}. This applies only if the components ${\bf f}^+\cdot\partial_{\bf x}\sigma$ and ${\bf f}^-\cdot\partial_{\bf x}\sigma$, of ${\bf f}^+$ and ${\bf f}^-$ normal to $\sigma=0$, have the same sign.

If, instead, the flow remains on the discontinuity, then \eref{inclusion} contains an element lying tangent to $\sigma=0$, given by
\begin{equation}\label{fs}
{\bf f}^{s}=\lambda^s{\bf f}^++(1-\lambda^s){\bf f}^-,\;\;\;\lambda^s:=\frac{{\bf f}^-\cdot\partial_{\bf x}\sigma}{({\bf f}^--{\bf f}^+)\cdot\partial_{\bf x}\sigma}.
\end{equation} 
Derived by Filippov \cite{f64,f88} (usually with \eref{inclusion} restricted to $\lambda\in[0,1]$), this prescribes the dynamics of the flow when sticking to the discontinuity. 
We need apply ${\bf f}^{s}$ only when both vector fields point towards, or both away from, the discontinuity, that is, when ${\bf f}^\pm\cdot\partial_{\bf x}\sigma$ have opposing signs, as shown at the shaded regions in \fref{fig:big2fold}.

Thus one obtains a flow that can hit the discontinuity and either cross it, concatenating solutions in the fields ${\bf f}^+$ and ${\bf f}^-$
, or stick to it, concatenating solutions in the fields ${\bf f}^\pm$ and ${\bf f}^{s}$
. The boundary between the two lies where ${\bf f}^+$ or ${\bf f}^-$ have turning points (tangencies) with respect to $\sigma=0$, where the flow {\it grazes} the discontinuity (see \fref{fig:big2fold}). 
The remainder of this letter deals with transitions between different behaviours at such boundaries. 

\begin{figure}[h!]\begin{center}\includegraphics[width=0.36\textwidth]{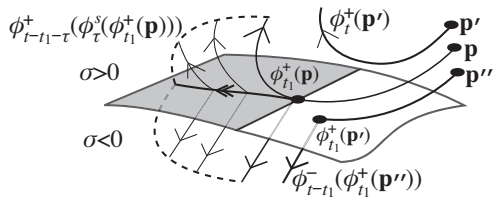}
\end{center}\vspace{-0.6cm}\caption{\sf Phase portrait of a grazing explosion. The flow through nearby points $\bf p$, ${\bf p}'$ and ${\bf p}''$, is found by concatenating the flows $\phi^+_t$, $\phi^-_t$, and $\phi^{s}_t$. An explosion takes place through $\bf p$.}\label{fig:grazing}\end{figure}

The flow consists of {\it all} solutions of \eref{inclusion} obtained by concatenating the flows of ${\bf f}^+$, ${\bf f}^-$, and ${\bf f}^{s}$. Consider what happens when the flow from a given point $\bf p$ grazes as in \fref{fig:grazing}. Under perturbation this grazing will not persist, meaning that if $\bf p$ is moved, the flow through $\bf p$ will generically either hit $\sigma=0$ transversally or miss it altogether. 

The case shown in \fref{fig:grazing} involves ${\bf f}^+$ curving away from $\sigma=0$, while ${\bf f}^-$ points away from it. The flow through $\bf p$ becomes multi-valued at grazing, in spite of which it can be precisely described as follows. 
Let $\phi^\alpha_t({\bf p})$ be the flow of ${\bf f}^\alpha$ through a point $\bf p$, so that
\begin{equation}
\frac{d\;}{dt}\phi^\alpha_t({\bf p})={\bf f}^\alpha\bb{\phi^\alpha_t({\bf p})},\quad \phi^\alpha_0({\bf p})={\bf p},
\end{equation}
where $\alpha$ denotes `$+$', `$-$', or `$s$'. 
The flow through $\bf p$ in $\sigma>0$ hits the discontinuity at time $t=t_1$ if $\sigma\bb{\phi^+_{t_1}({\bf p})}=0$, and grazes if ${\bf f}^+\cdot\partial_{\bf x}\sigma\bb{\phi^+_{t_1}({\bf p})}=0$. For $t<t_1$ the flow is single-valued. For $t>t_1$ the state undergoes an explosion of possible values, given by 
\begin{equation}\label{unfoldset}
{\bf x}\in\cc{{\bf x}_\tau=\phi^\pm_{t-t_1-\tau}\bb{\phi^{s}_{\tau}\bb{\phi^+_{t_1}({\bf p})}}:\tau\in[0,t-t_1]}.
\end{equation}
Each orbit in the flow sticks for a time $\tau\in[0,t-t_1]$ after grazing, then leaves $\sigma=0$ following either $\phi^+_t$ or $\phi^-_t$. If $\bf p$ is slightly shifted (to ${\bf p}'$ or ${\bf p}''$), for instance because a parameter in the system changes, the flow switches discontinuously between ${\bf x}=\phi^+_{t}({\bf p}')$ and ${\bf x}=\phi^-_{t-t_1}\bb{\phi^+_{t_1}({\bf p}'')}$, in the latter case hitting the discontinuity transversally at time $t=t_1$. This is classified as a form of bifurcation in \cite{jh09}, neglecting the explosion \eref{unfoldset} itself. 

The significance of such an event is seen if either ${\bf p}'$ or ${\bf p}''$ lies on, for example, a limit cycle in the flow. If varying a parameter causes grazing in the manner of \fref{fig:grazing}, the simple local geometry reveals that the cycle must suffer a sudden explosion, followed by abrupt disappearance of both the explosion and the cycle itself. The derivation of this event, forthcoming in \cite{jh09}, was promptly followed by its observation as the mechanism for sudden onset of thermal oscillations in a superconductor \cite{js09}, for which only preliminary studies of a lower dimensional approximate model have previously been made. Below, the grazing event is identified and simulated in the full model. 

\begin{figure}[b!]\centering\includegraphics[width=0.47\textwidth]{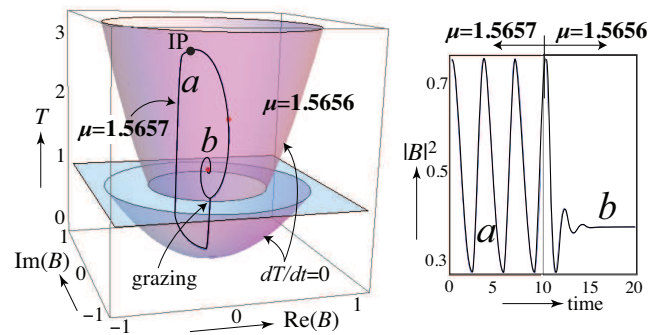}
\vspace{-0.3cm}\caption{\sf Effect of grazing explosion in a superconductor: simulation of \eref{resonator} for $\varepsilon=100$, $s^+=3.891$, $s^-=1.297$, $\Lambda^-=-0.2+i$, $\Lambda^+=-0.5+i\mu$. For two nearby values of $\mu$ (see picture), an orbit through the initial point (IP at $B=0.8-0.4i$, $T=3$) oscillates around the cycle $a$, crossing between $T>1$ and $T<1$, or finds a steady state $b$ in $T>1$. The plane where $T=1$ and parabolae where $d T/dt=0$ are shown. In the right figure, the time trace of the power $|B|^2$ is shown, with the small change in $\mu$ introduced at time $t=10$. }\label{fig:super}\end{figure}

{\it An experimentally motivated example: explosion in a superconducting resonator.  } 
The device of interest is a stripline resonator designed for high sensitivity measurements of quantum phenomena, such as the Casimir effect \cite{bs08}. 
With $i=\sqrt{-1}$, and $\varepsilon$ a small parameter, the system
\begin{equation}\label{resonator}
\begin{array}{rcl}\frac{d\;}{dt} B&=&\Lambda(T) B-{\rm i},\\\varepsilon\frac{d\;}{dt}  T&=&s(T)|B|^2-T,
\end{array}
\end{equation}
describes the complex current amplitude, $B$, in a ring of niobium nitride at temperature $T$, which is superconducting when sufficiently cold. 
All quantities are non-dimensionalised so that the ring is superconducting for $T<1$ and normal conducting for $T>1$, causing $\Lambda$ and $s$ to jump between constant values $\Lambda^\pm$ and $s^\pm$,
\begin{equation}\label{respws}
(\Lambda(T),s(T))=\left\{\begin{array}{lll}(\Lambda^+,s^+)&\mbox{if}&T>1,\\(\Lambda^-,s^-)&\mbox{if}&T<1.\end{array}\right.
\end{equation}
Grazing occurs when \eref{resonator} lies tangent to the discontinuity surface $T=1$, hence where $\frac{d\;}{dt} T=0$. \Fref{fig:super} shows a limit cycle on which these conditions are satisfied, along with a small perturbation that destroys the cycle in the manner of \fref{fig:grazing}, after which the system settles rapidly onto a pre-exisiting steady state. 
The effect fits qualitatively with experimental observations (see \cite{bs08,js09} and references therein). Further work remains, to measure the parameter values at which grazing occurs experimentally, and to study the way in which, near grazing, noise causes repeated jumps between the modes $a$ and $b$ in \fref{fig:super}. 

It follows from the geometry in \fref{fig:grazing} (see \cite{jh09}) that grazing explosions are codimension one phenomena, meaning that generically they may be observed as a parameter varies. However, this also implies that the many outcomes possible during the explosion (shown in \fref{fig:grazing} but omitted between $a$ and $b$ in \fref{fig:super}) are not readily observed in practice. We now show that there are scenarios where the full explosion is observable even under perturbation, with striking consequences. 
%
%
In more than two dimensions there may exist points where both ${\bf f}^+$ or ${\bf f}^-$ are tangent to the discontinuity simultaneously, 
which allows non-determinism to affect the flow more robustly than in the grazing explosion.  

\begin{figure}[!h]\centering\includegraphics[width=0.26\textwidth]{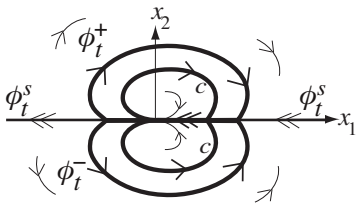}
\vspace{-0.3cm}\caption{\sf Non-deterministic chaos in the system \eref{dbfold}. For any two points outside the arcs $c$, there exists an orbit from the origin passing through them. We can view this as extreme sensitivity to initial conditions caused by non-determinism at the double fold. }\label{fig:vischaos}\end{figure}

Consider first the scenario of a planar system that is symmetric about a discontinuity, given by \eref{pws} with
\begin{equation}\label{dbfold}
{\bf f}^\pm=\left(\begin{array}{c}-1\pm x_2\\x_2\mp x_1\end{array}\right),\quad\sigma(x_1,x_2)=x_2,
\end{equation}
for which the points $(x_1,x_2)=(1,\pm1)$ are unstable foci. Then, on $\sigma=0$, ${\bf f}^\pm\cdot\partial_{\bf x}\sigma=0$ at the origin, and ${\bf f}^s=(-1,0)$. 
As a simple sketch shows (\fref{fig:vischaos}), every point in the flow reaches the origin in forward time, and may do so via sticking. Thereafter the flow consists of sequences of arcs $\phi^\pm_t$ outside the discontinuity, and sticking $\phi^{s}_\tau$ for indeterminate times $\tau$. Thus the system makes recurrent visits to the origin, interspersed by excursions into $\sigma\neq0$ with unpredictable durations and trajectories. 

The model above is artificial because of the imposed symmetry, but such scenarios do occur generically in higher dimensions, at transverse intersections between sets ${\bf f}^\pm\cdot\partial_{\bf x}\sigma=0$ where ${\bf f}^+$ or ${\bf f}^-$ are quadratically tangent to the discontinuity (as in \fref{fig:big2fold}). At such points, two different regions of sticking meet, and by concatenating the flows of ${\bf f}^\pm$ and ${\bf f}^{s}$, one obtains a flow that can pass through the double tangency from attracting to repelling regions of the discontinuity surface, similar to \fref{fig:vischaos}. 

\begin{figure}[h!]\begin{center}\includegraphics[width=0.47\textwidth]{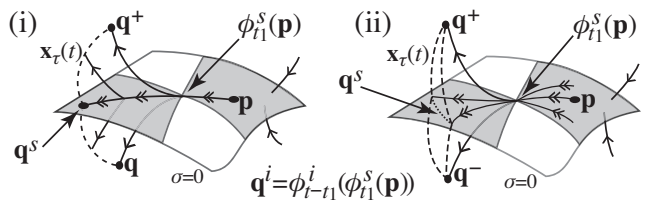}
\end{center}\vspace{-0.5cm}\caption{\sf Explosion at a double tangency: the flow $\phi^{s}_t({\bf p})$ hits the singularity at $t=t_1$. The subsequent multi-valued flow is given by \eref{un2foldset}, either sticking along $\phi^{s}_t$ or escaping along $\phi^\pm_t$. In (i) only one, and in (ii) many, orbits traverse a single point between regions where the discontinuity changes from attracting to repelling. Some particular points ${\bf q}^i$ are labelled.}\label{fig:un2fold}\end{figure}

\Fref{fig:un2fold} shows scenarios, generic in three or more dimensions \cite{f88}, where one (i) or many (ii) orbits pass between sticking regions via a double tangency. In the figure, the flow though a point $\bf p$ is single-valued until it hits a double tangency at time $t_1$, after which the state explodes into a family of orbits that stick for a time $\tau\in[0,t-t_1]$, then follow ${\bf f}^+$ or ${\bf f}^-$ until time $t$, expressed as 
\begin{equation}\label{un2foldset}
{\bf x}\in\cc{\;{\bf x}_\tau=\phi^\pm_{t-t_1-\tau}\bb{\phi^{s}_{\tau}\bb{\phi^{s}_{t_1}({\bf p})}}:\tau\in[0,t-t_1]}.
\end{equation}

The scenario in \fref{fig:un2fold}(ii) makes possible an extreme manifestation of non-determinism, similar to \fref{fig:vischaos}. If some global mechanism exists that returns the flow through ${\bf q}^\pm$ to a neighbourhood of $\bf p$, then a set is generated in which the flow returns recurrently to the double tangency, yet each excursion has an unpredictable duration and trajectory, constituting non-deterministic chaos. A formal definition of this requires augmenting the definition of {\it deterministic} chaos, by extending to multi-valued flows the idea of sensitive dependence on initial conditions. Such an extension is given in \cite{cj10}, where non-deterministic chaos arises as a relatively rare event near a novel bifurcation.
We conclude with a potentially more common example that exhibits the phenomenon. 

\begin{figure}[h!]\centering\includegraphics[width=0.47\textwidth]{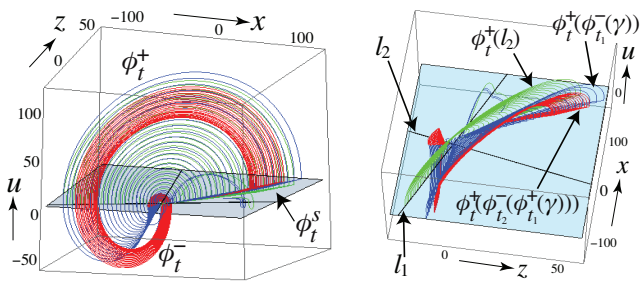}
\vspace{-0.3cm}\caption{\sf Non-deterministic chaos in a simulation of \eref{mech}. The two views show the discontinuity surface $u=0$ and tangency lines $l_{1,2}$. The flows in $u>0$, $u<0$, and the sticking flow on $u=0$, are indicated by $\phi^+_t$, $\phi^-_t$, $\phi^{s}_t$. Several orbits are shown on the boundaries of the region of non-deterministic chaos.}\label{fig:mechns}\end{figure}

{\it A toy model: explosion in a mechanical oscillator.  } 
Consider an object of unit mass, whose displacement $x$ satisfies a Newtonian force law $\ddot x=(\dot x-v)b-x+g(\dot x,t)$, where $\dot x=dx/dt$. This includes a spring force $-x$, a negative damping proportional to the speed relative to some reference $v$, plus an additional forcing $g$. For $\dot x<v$ let $g$ grow linearly in time, say as $g=r_1t$. For $\dot x>v$ let $g$ have speed-dependent dynamics, setting $g=r_2z$ where $\dot z=a+(\dot x-v)c$. For this abstract illustration we let $a=-1.3$, $b=0.1$, $c=0.2$, $v=-1$, $r_1=12$, $r_2=-1$. 
Letting $u=\dot x-v$, we obtain a first order system
\begin{equation}\label{mech}
\begin{array}{rcl}\frac{d\;}{dt}  x&=&u+v\\\frac{d\;}{dt}  u&=&-x+bu+r_1z+(r_2-r_1)z\unitstep{u},\\\frac{d\;}{dt}  z&=&1+(a-1+cu)\unitstep{u},
\end{array}
\end{equation}
where $H(u)=1$ for $u>0$ and $H(u)=0$ for $u<0$. The flows are tangent to the discontinuity surface, $u=0$, along the lines $l_{1,2}$ where $x=r_{1,2}z$ (labelled in \fref{fig:mechns}), which intersect at the origin. Part of the flow through the origin detaches from $u=0$, winds around through both $u>0$ and $u<0$, returns to $u=0$, then sticks and returns to the origin, local to which the phase portrait resembles \fref{fig:un2fold}(ii). This generates a non-deterministic chaotic set, on which the entire flow meets the double tangency both in forward and backward time, and which is highly attracting with respect to the surrounding flow. 

The term `non-deterministic chaos' is not intended to suggest that smoothing the discontinuity will yield deterministic chaos. It is a simple exercise to simulate \eref{mech} with the step $\unitstep u$ replaced by a sigmoid function, which yeilds a system where the flow explores the set shown in \fref{fig:mechns}, but is enormously sensitive to factors such as numerical precision and the method of smoothing. Initial study of this and other examples suggests that smoothing can lead to chaos, canards, mixed mode oscillations, or simple limit cycles. Further study is warranted. 


Discontinuous models are increasingly common in mechanical, electronic, and biological applications, often approximating sigmoidal changes for which a precise description is lacking. A general understanding of discontinuity induced phenomena is slowly emerging. In many cases, dynamics is robust to the precise form of the sigmoid/discontinuous jump. Explosions, on the other hand, will certainly be sensitive to the jump's precise form. The multi-valued flow evidently plays a role in organising non-determinism in the discontinuous limit. 

The genericity of the vector fields underlying \fref{fig:grazing} and \fref{fig:un2fold} has been understood since the seminal work of Filippov \cite{f88}, but their non-deterministic consequences have received little attention. Although illustrated here in three dimensions, these derive from normal form vector fields generic in any higher dimension \cite{jh09}. When dimensions are added, the local geometry responsible for non-determinism remains, the set-valued flows retain the same dimension, and the phenomena remain of codimension one (\Frefs{fig:grazing} \& \ref{fig:un2fold}(i)) or codimension zero (\Fref{fig:un2fold}(ii)), provided that discontinuities of the local form (\ref{pws}) exist. Explosions, as introduced here, reveal extreme unpredictability caused not by external noise or defective modeling, but by inherent local geometry. With this fact having been largely overlooked, it is not surprising that examples beyond those above remain to be uncovered. 

It should be noted that the flows considered here are {\it continuous}, with a discontinuous time derivative. Non-uniqueness in forward time is not uncommon in systems involving hybrids of flows with maps, such as impacts or finite state resets, where the flow itself is discontinuous. In the present study it is the continuity of the flow, coupled with loss of determinism, that generates explosions associated with tangency to a discontinuity.

A discontinuity in the velocity field of a continuous flow allows multiple histories or outcomes from a single point. 
Regarding the latter, it might be argued that any physical application must contain information that restores determinism. Such information is not, however, part of the model in the discontinuous limit. 
Nor is it guaranteed that one has practical access to such information that is not overwhelmed by noise or uncertainty. The discontinuous model provides a geometric description in the limit where infinite sensitivity to initial conditions constitutes the breakdown of determinism.


\bibliography{../../grazcat}
\bibliographystyle{plain}

\end{document}